\documentclass[11pt,english]{article}

\usepackage{anysize}
\marginsize{2cm}{2cm}{2.5cm}{2.4cm}

\usepackage[english]{babel}
\usepackage{amsfonts}
\usepackage{latexsym}
\usepackage{amsmath}
\usepackage{mathrsfs}
\usepackage{amssymb}
 \usepackage{graphpap}
 \usepackage{graphicx}
 \usepackage{epsfig}

\let\phi=\varphi
\newcommand{\Z}{{\mathbb Z}}
\newcommand{\R}{{\mathbb R}}
\newcommand{\N}{{\mathbb N}}
\newcommand{\eps}{\varepsilon}
\newcommand{\IP}{{\mathbb P}}
\newcommand{\IE}{{\mathbb E}}

\newcommand{\Po}{{\mathtt P}_{\omega}}

\newcommand{\Eo}{{\mathtt E}_{\omega}}

\newcommand{\sig}{\sigma}
\newcommand{\alf}{\alpha}

\newcommand{\qed}{\hfill$\Box$\par\medskip\par\relax}

\newcommand{\argmax}{\mathop{\mathrm{arg\,max}}}

\let\phi=\varphi

\newcommand{\lf}{\lfloor}
\newcommand{\rf}{\rfloor}
\newcommand{\lc}{\lceil}
\newcommand{\rc}{\rceil}

\newcommand{\om}{{\omega}}

\def\P{\mathbf{P}}

\newtheorem{theo}{Theorem}[section]
\newtheorem{lm}{Lemma}[section]

\newtheorem{prop}{Proposition}[section]
\newtheorem{cor}{Corollary}[section]

\allowdisplaybreaks

\title{Localization for a random walk in slowly decreasing random potential}
\author{Christophe Gallesco$^{~1}$ \and Serguei
Popov$^{~2}$ \and Gunter M.\ Sch\"utz$^{~3}$}

\begin{document}

\bibliographystyle{plain}

\maketitle
{\footnotesize 

\noindent $^{1,2}$Department of Statistics, 
Institute of Mathematics, Statistics and Scientific Computation,
University of Campinas--UNICAMP, 
rua S\'ergio Buarque de Holanda 651, 13083--859, Campinas SP,
Brazil\\
\noindent e-mail: \texttt{gallesco@ime.unicamp.br}\\
\noindent e-mail: \texttt{popov@ime.unicamp.br}; 
URL: \texttt{http://www.ime.unicamp.br/$\sim$popov}\\
\noindent $^{~3}$ Theoretical Soft Matter and Biophysics, Institute of Complex Systems, Forschungszentrum J\"ulich, \\
52425 J\"ulich, Germany\\
\noindent e-mail: \texttt{g.schuetz@fz-juelich.de}
}

\maketitle

\begin{abstract}
We consider a continuous time random walk $X$ in random environment on $\Z^+$ such that its potential can be approximated by the function $V: \R^+\to \R$ given by $V(x)=\sig W(x) -\frac{b}{1-\alf}x^{1-\alf}$ 
where $\sig W$ a Brownian motion with diffusion coefficient $\sig>0$ and parameters $b$, $\alf$ are 
such that $b>0$ and $0<\alf<1/2$.
We show that $\P$-a.s.\ (where $\P$ is the averaged law) 
$\lim_{t\to \infty} 
\frac{X_t}{(C^*(\ln\ln t)^{-1}\ln t)^{\frac{1}{\alf}}}=1$ with
$C^*=\frac{2\alf b}{\sig^2(1-2\alf)}$.
In fact, we prove that by showing that there is a trap located
around $(C^*(\ln\ln t)^{-1}\ln t)^{\frac{1}{\alf}}$ (with 
corrections of smaller order) where
the particle typically stays 
up to time~$t$. This is
in sharp contrast to what happens in the ``pure'' Sinai's regime,
where the location of this trap is random on the scale $\ln^2 t$.
\\[.3cm] 
\\[.3cm] 
\end{abstract}

\section{Introduction and results}
\label{s_intro}
Suppose that $\omega=(\omega_x)_{x\geq 1}$ is a sequence of a i.i.d.\ random variables.
Fix $b>0$ and $\alf \in (0,\frac{1}{2})$ and let us define the sequence $(q_y)_{y\geq 0}$ such that 
$q_0=0$ and $q_y=\frac{\exp(\om_y-by^{-\alf})}{1+\exp(\om_y-by^{-\alf})}$ for $y\geq 1$. For each 
realization of~$\om$, we consider the continuous time random walk~$X$ 
on~$\Z^+$ with transition probabilities given by
\begin{align*}
\Po[X_{t+h}&=y+1\mid X_t=y]=(1-q_y)h+o(h),\nonumber\\
\Po[X_{t+h}&=y-1\mid X_t=y]=q_yh+o(h),\phantom{*****} \mbox{if $y\geq 1$},
\end{align*}
as $h\to 0$. We will denote by
$\IP, \IE$ the probability and expectation with respect to
$\omega$, and by  $\Po$, $\Eo$  the (so-called ``quenched")
probability and expectation for the random walk in the fixed
environment $\omega$. We will use the notation $\Po^{x}$ for the 
quenched law of $X$ starting from~$x$. Nevertheless, for the sake of brevity, 
we will omit the superscript $x$ whenever $x=0$.
We make the following assumption :

\medskip
\noindent
\textbf{Condition S.} We have
\[
\IE[\omega_1]= 0,\quad
\sigma^2 := \IE[\omega_1^2]
\in (0, +\infty).
\]

\medskip

The vanishing expectation of  $\omega_1$ means that the random walk has a drift which is asymptotically decaying,
which is the case of interest to be studied here.
For technical reasons we also assume that the following condition holds:

\medskip
\noindent
\textbf{Condition K.} There exists a $\theta_0>0$ such that $\IE[e^{\theta\om_1}]<\infty$ for all $|\theta|<\theta_0$.

\medskip

The choice of the rates $q_y$ has the interpretation of a random walk in a power law potential with amplitude $b$
on which a Sinai-type random potential is superimposed. Indeed, in the case $b=0$, Condition~S corresponds to 
Sinai's regime~\cite{Sinai} (after stating our main result, we will compare
it with what happens in ``pure'' Sinai's regime). 
Random walks in an asymptotically decaying power-law potential 
play an important role 
in a number of applications in physics. As a very well-studied example we mention the condensation 
transition in the zero-range process where the grand-canonical stationary distribution on a single site 
is that of a random walk in a power-law (or
logarithmic) potential \cite{ZRP1,ZRP2,ZRP3,ZRP4}. 
For $0<\alf<1$ and $b<0$ there exists a finite critical particle density above which the 
grand-canonical
stationary distribution does not exist. 
Then, in a canonical ensemble with fixed total particle number 
such that the total density exceeds the critical value, a macroscopic
number of particles ``condenses'' on a single site. 
The same is true for $\alf=1$ and $b\leq -2$, a case
of particular importance e.g.\ in DNA denaturation where by a mapping to the dynamics of unzipped DNA
strands the presence or absence of a 
condensation transition indicates whether the DNA denaturation transition is of first or second 
order \cite{Kafr02,Ambj06}. It is then natural to study the effect of quenched disorder which is usually 
modelled by a random potential of the 
type defined above. It turns out that the condensation transition persists only in the range $0<\alf<1/2$ \cite{CGS},
which appears to be related to the smoothening of depinning transitions for directed polymers with quenched disorder
of which the DNA denaturation transition is an example \cite{GT1,GT2}.

Directly from the viewpoint of random walks in random environments 
the presence of quenched disorder 
in an asymptotically decaying power-law potential
has been studied in detail in \cite{Menswade1,Menswade2} in a discrete time setting.
The presence of a 
condensation transition corresponds to ergodicity of the random walk. Going beyond stationary properties,
these authors relate the position of the random walk to some expected hitting times
to obtain a series of interesting results on the speed of the random walk starting from the origin.
In this respect the transient case is of particular interest. For $b>0$ and $\alf \geq 1/2$ the scenario is not very much different from
the case of pure Sinai-disorder (no power law potential). Roughly speaking, the displacement of the random walk
from the origin grows to leading order in time $t$ as $(\ln t)^2$, independent of $\alf$. On the other hand, for $b>0$ and 
$0 < \alf < 1/2$ it was proved \cite{Menswade2} that for 
a.e.\ random environment~$\omega$ one has a.s.\
$(\ln \ln t)^{-1/\alf - \epsilon} < \eta_t(\omega)/(\ln t)^{1/\alf} < (\ln \ln t)^{2/\alf + \epsilon} $
for all but finitely many $t$.

The approach used here allows us to go further. The main result of this paper is:
\begin{theo}
\label{Theo}
Under Conditions~S and K, we have for $\IP$-almost all realizations of $\omega$,
\begin{equation*}
 \lim_{t\to \infty}\frac{X_t}{(C^*(\ln\ln t)^{-1}\ln t)^{\frac{1}{\alf}}}=1,\phantom{***}\mbox{$\Po$-a.s.,}
\end{equation*}
with $C^*= \frac{2\alf b}{\sig^2(1-2\alf)}$.
\end{theo}

Observe that we define the model in a continuous-time 
setting rather than in discrete time. This brings about a (very) slight technical complication, 
but is better motivated from a physics perspective.

Let us comment now on the relationship of our work with the classical
model of one-dimensional RWRE in i.i.d.\ environement (see e.g.\ \cite{Zeitouni}).
As often happens with theorems of this kind, the proof of Theorem~\ref{Theo}
is obtained by showing that the particle will eventually find a \textit{trap}
(i.e., a piece of the environment with ``drift inside''), and then
 stay there up to time~$t$. It it well-known that, for the RWRE in
Sinai's regime, the location of this trap (scaled by $\ln^2 t$) is a random variable. However, it is interesting to observe that (as one can see from the proof of Theorem~\ref{Theo}) adding the power-law perturbation to the 
Sinai's potential changes the situation: the position of \textit{the trap}
becomes ``less random'' (there are still fluctuations, of course, but they are
of smaller order).

As an aside, we mention that with
%\textit{Remark.} Define 
$s(t):=(C^*(\ln\ln t)^{-1}\ln t)^{\frac{1}{\alf}}$ we can also 
deduce from the proof of Theorem~\ref{Theo}, the following upper bounds for some particular hitting times of $X$.
For $\eps\in (0,1)$, let $\tau_{(1-\eps)s(t)}$ be the first hitting time of the point $\lf(1-\eps)s(t)\rf$ by the random walk $X$. 
Then, for all $\eps>0$ there exists $\delta>0$ such that $\IP$-a.s., 
\[
\Po[\tau_{(1-\eps)s(t)}>t]\leq \exp\{-t^{\frac{\delta}{2}}\}
\]
for all $t$ large enough (see equation~(\ref{FUG1})).
\medskip

In the next section, we introduce some notations and recall some auxiliary facts which 
are necessary for the proof of Theorem~\ref{Theo}. In section \ref{Technical}, we prove various 
technical lemmas about the asymptotic behavior of the environment. Finally, in section \ref{secTheo}, 
we give the proof of Theorem~\ref{Theo}.

\section{Notations and auxiliary facts}

Given a realization of~$\omega$, define the potential function
for $x\in\mathbb{R^+}$, by
\[
U(x) :=\sum_{y=1}^{\lf x\rf}\ln\frac{q_y}{1-q_y}=\sum_{y=1}^{\lf x\rf}(\om_y-by^{-\alf})
\]
where $\lf x\rf$ is the integer part of $x$ and $\sum_{y=1}^{\lf x\rf}:=0$ if $x<1$. The behavior of  $U$ is of crucial importance for the analysis of the asymptotic properties of the random walk $X$ (cf.\ Propositions \ref{Confine} and \ref{escape} below). 

\begin{figure}[!htb]
\begin{center}
\includegraphics[scale= 0.6]{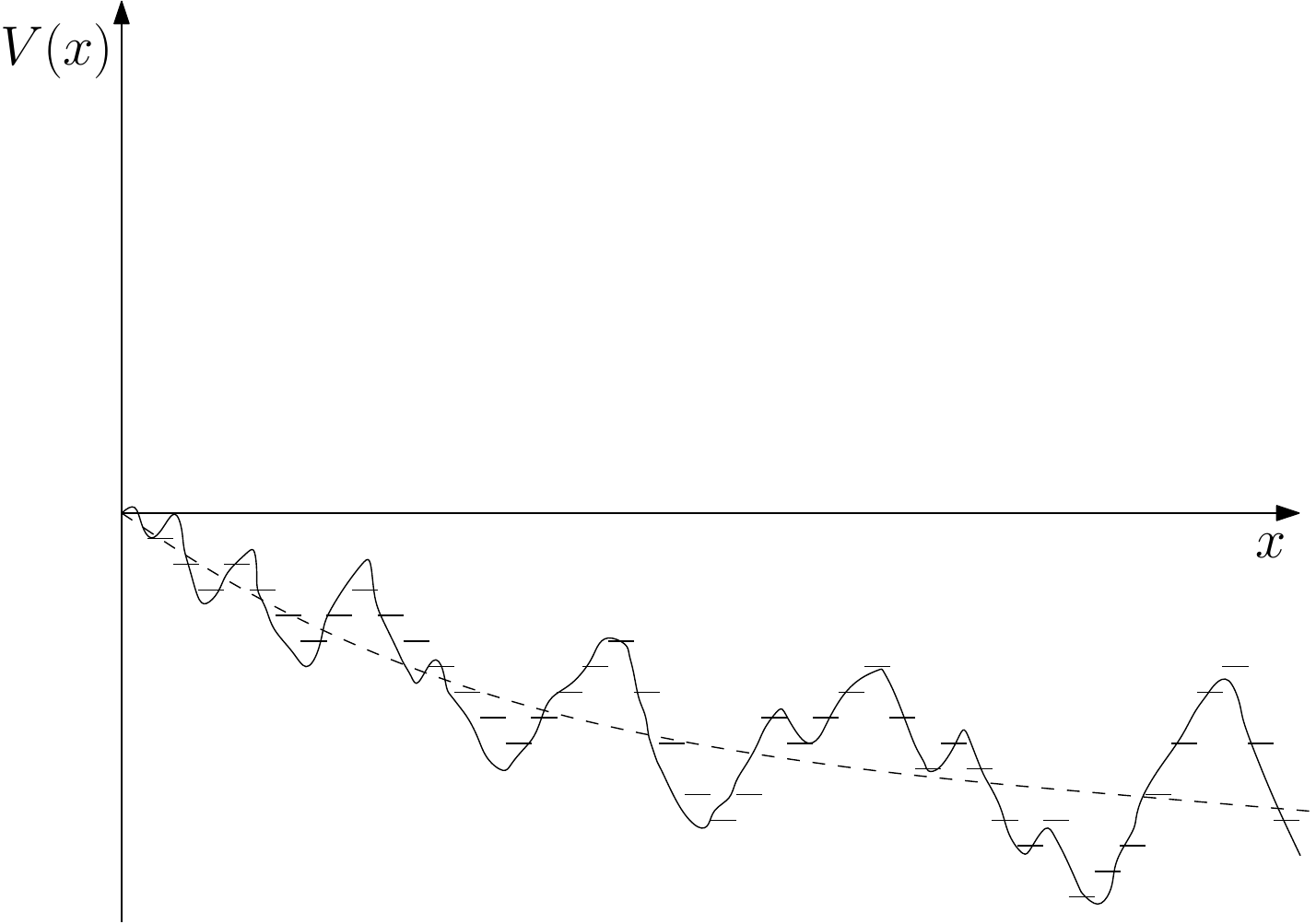}
\caption{Approximation of potential $U$ by $V$.}
\label{fig1}
\end{center}
\end{figure}

Conditions S and K will allow us to couple the potential $U$ to Brownian motion with power law drift, simplifying much the proof of limit properties of the random walk $X$.
Indeed, by the well-known Koml{\'o}s-Major-Tusn{\'a}dy strong approximation
theorem (cf.\ Theorem 1 of \cite{KMT}), there exists (possibly in an enlarged probability space) a coupling for $\om$ and a standard Brownian motion $W$, such that 
\begin{equation}
\label{KMT}
\IP\Big[\limsup_{n\to \infty}\frac{\max_{1\leq m\leq n}|\sum_{i=1}^m\omega_i-\sig W(m)|}{\ln n}\leq \hat{K}\Big]=1
\end{equation}
for some finite constant $\hat{K}>0$. 
A useful consequence of (\ref{KMT}) is that if $x$ is not too far away from the origin, then $\sum_{i=1}^{\lf x\rf}\om_i$ and $\sig W(x)$ are rather close for the vast majority of environments. Hence, it
is convenient to introduce the following set of ``good'' environments and 
to restrict our forthcoming computations to this set.
Fix $M>\frac{1}{\alf}$ and for any~$t>e$, let
\begin{equation}
\label{approx1}
\Gamma(t): = \Big\{\omega : \Big|\sum_{i=1}^{\lf x\rf} \om_i-\sig W(x)\Big|\leq K\ln\ln t\;, \; x\in [0,\ln^{M}t]\Big\}.
\end{equation}
By (\ref{KMT}) and properties of the modulus of continuity of Brownian motion, we can choose $K \in (0, \infty)$ in such a way that for $\IP$-almost all~$\omega$, it holds that $\omega\in\Gamma(t)$
for all~$t$ large enough (cf. e.g. \cite{CP} or \cite{Galles} where this fact was used).
On the other hand, using the fact that there exists a finite constant $C>0$ such that for all $x\geq 1$, 
\begin{equation}
\label{approx2}
\Big|\sum_{i=1}^{\lf x\rf}i^{-\alf}-\int_1^xu^{-\alf}du\Big|\leq C,
\end{equation}
we can define a new potential function $V$ by 
\[V(x):=\sig W(x)-\frac{b}{1-\alf}x^{1-\alf}\] 
for all $x\in \R^+$ and using~(\ref{approx1}) 
and~(\ref{approx2}), we have that there exists a finite $K_1>0$ such that for all $t> e$ and $\om \in \Gamma(t)$, $\max_{x\leq \ln^M t}|V(x)-U(x)|\leq K_1\ln\ln t$. Observe that $V$ is a Brownian motion with a power law drift. For convenience, from now on, we will work with potential $V$ instead of $U$ (see Fig.\ \ref{fig1}).

\begin{figure}[!htb]
\begin{center}
\includegraphics[scale= 0.6]{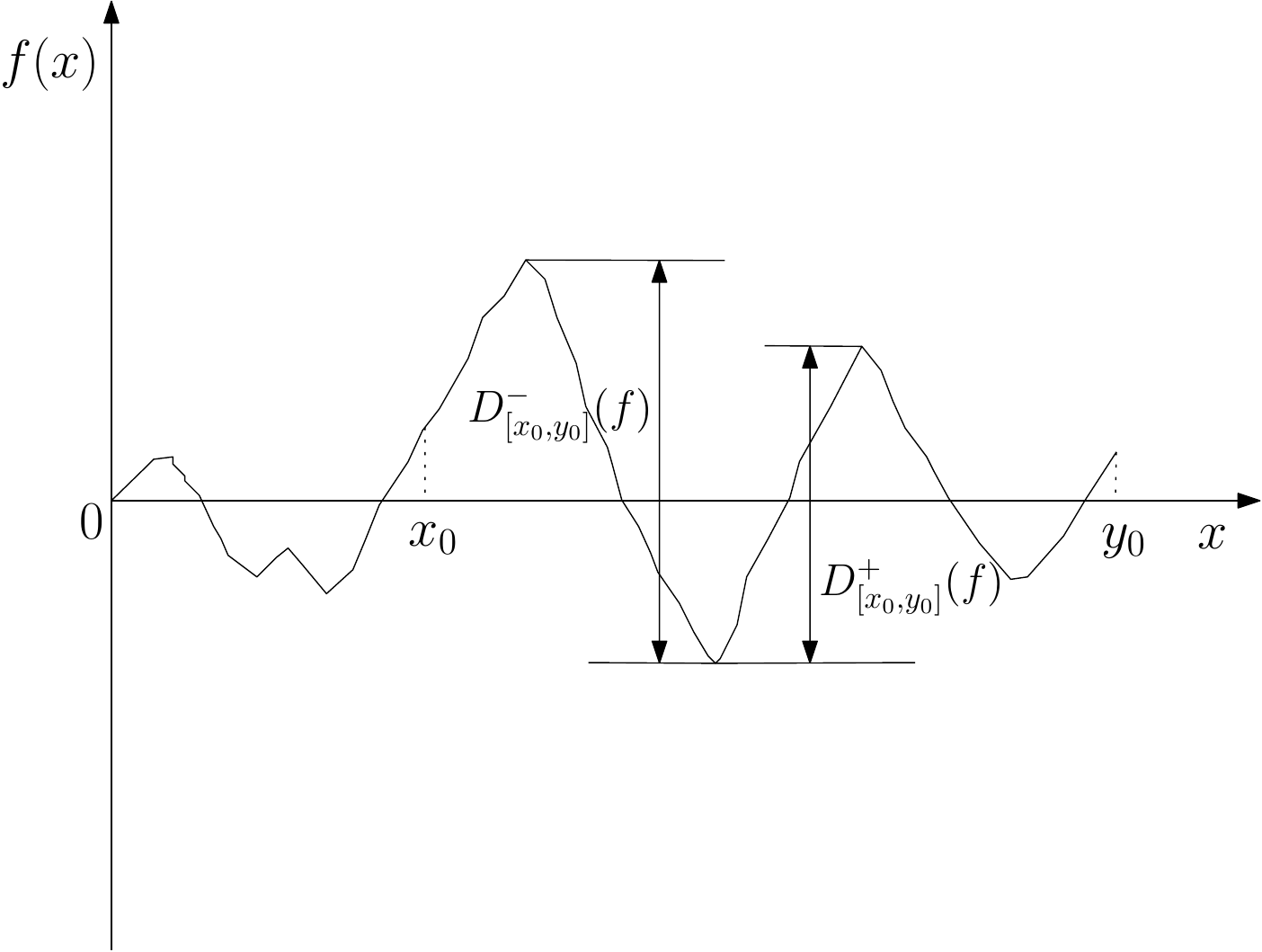}
\caption{On the definitions of $D^+_{[x_0,y_0]}(f)$ and $D^-_{[x_0,y_0]}(f)$.}
\label{fig1b}
\end{center}
\end{figure}

\medskip

For a function $f:\R^+\to \R$ and $x_0<y_0$, let $D^+_{[x_0,y_0]}(f):=\sup_{u\in [x_0,y_0]}(f(u)-\inf_{v\in[x_0,u]}f(v))$ and $D^-_{[x_0,y_0]}(f):=\sup_{u\in [x_0,y_0]}(f(u)-\inf_{v\in[u,y_0]}f(v))$ be respectively the maximum draw-up and draw-down of the function $f$ on the interval $[x_0,y_0]$ (see Fig.~\ref{fig1b}). As we will see in the proof of 
Theorem~\ref{Theo}, these functionals applied to the potential $V$ are important quantities in order to determine the limiting behavior of the random walk $X$. 
The distribution of $D_{[x_0,y_0]}^+$ is not known for a Brownian motion with a power law drift. Fortunately, in our case, we can locally approximate the power law drift by a linear one. 
It happens that for fixed intervals~$I$ the law of $D_I^+$ is known for a Brownian motion 
with linear drift (cf.\ (1) 
in~\cite{MAPAM}) but  in this reference, it is given under the form of an alternating 
series which is not easy to handle. If, instead of considering deterministic intervals~$I$ 
we consider intervals of 
size given by an exponential random variable independent of $W$ then the law of $D_I^+$ 
becomes much simpler and is more useful for our purposes. 
\medskip

We now recall the following result which can be found in \cite{Salminen}: 

\begin{prop}
\label{Prop0}
Let $T$ be a random variable with exponential distribution of mean $\mu$ and $W^{(\sig,\nu)}$ a Brownian motion with diffusion coefficient $\sig$ and linear drift $\nu$, that is, $W^{(\sig,\nu)}(t)=\sig W(t)+\nu t$ where $W$ is a standard Brownian motion. Assume that $T$ is independent of $W$.
Then, 
\begin{equation*}
P\Big[D_{[0,T]}^+(W^{(\sig,\nu)})>a\Big]=\frac{\exp(\nu a\sig^{-2})}{\cosh(a\sig^{-1}\sqrt{2\mu^{-1}+\nu^2\sig^{-2}})+\frac{\nu\sig^{-1}}{\sqrt{2\mu^{-1}+\nu^2\sig^{-2}}}\sinh(a\sig^{-1}\sqrt{2\mu^{-1}+\nu^2\sig^{-2}})}
\end{equation*}
for all $a\geq0$.
\end{prop}

It is then not difficult to establish the following
\begin{cor}
\label{Corro1}
Suppose that $\nu<0$ and that $a$, $\nu$ and $\mu$ are functions the real variable $t>0$. 
If $a|\nu|\to \infty$, $\nu^2\mu \to \infty$ and $a(\mu |\nu|)^{-1}\to 0$ as $t\to \infty$, 
then
\begin{equation*}
P\Big[D_{[0,T]}^+(W^{(\sig, \nu)})>a\Big]= \frac{1}{1+\frac{\sig^2}{2\nu^2\mu}\exp(\frac{2|\nu|a}{\sig^2})}(1+o(1))
\end{equation*}
as $t\to \infty$.
\end{cor}

For all $A\subset \Z^+$ we define $\tau_A:=\inf\{t>0: X_t\in A\}$
the first hitting time of $A$ for the random walk~$X$. When $A=\{x\}$, $x\in \Z^+$, we simply write $\tau_x$ instead of $\tau_{\{x\}}$.\\
Let $I=[a,b]$ with $0\leq a<b<\infty$ be a finite interval of $\Z^+$ and let $H(I):=D^+_I(U)\wedge D^-_I(U)$ and $\tilde{M}:=D^+_I(U)\vee D^-_I(U)$. We will need the following upper bound on the probability of confinement which comes from the proof of Proposition 4.1 of \cite{PGF}:
\begin{prop}
\label{Confine}
There exists a positive constant $K_2$ such that, $\IP$-a.s., for any finite interval $I=[a,b]$ and any point $x$ such that $a<x<b$,
\begin{equation*}
\Po^x[\tau_{\{a,b\}}\geq t] \leq \exp \Big\{-\frac{t}{K_2(b-a)^3(b-a+\tilde{M})e^{H(I)}}\Big\}
\end{equation*}
for all $t>K_2(b-a)^3(b-a+\tilde{M})e^{H(I)}$.
\end{prop}

For the random walk $X$, we will
 eventually need to estimate the probability of escaping to one 
specific direction. In Proposition \ref{escape}, as an example, we just state the result for the probability of escaping to the right. Nevertheless, in section \ref{secTheo}, we will use this estimate in both directions. 
We define a reversible measure $\pi$ by $\pi(0):=1$ and $\pi(x):= e^{-U(x)}+e^{-U(x-1)}$ for $x\geq 1$ (observe that $\pi(x)(1-q_x)=q_{x+1}\pi(x+1)$ for all $x\in \Z^+$). 
For any finite interval $I$ of $\Z^+$, we define $h_I:=\argmax_{x\in I}U(x)$. 
We will use the following estimate (see e.g.\  the proof of \ Proposition 4.2 in \cite{PGF}):
\begin{prop}
\label{escape} 
There exists a positive constant $K_3$ such that, $\IP$-a.s., for any finite interval $I=[a,b]$ of $\Z^+$ we have
\begin{equation*}
\Po^a[\tau_b<t] \leq K_3t\frac{\pi(h_I)}{\pi(a)}
\end{equation*}
 for all $t>1$.
\end{prop}
Using the above expression of the reversible measure $\pi$, we have 
\begin{equation*}
\frac{\pi(h_I)}{\pi(a)}\leq e^{-U(h_I)+U(a)}\Big(1+e^{U(h_I)-U(h_I-1)}\Big).
\end{equation*}
If $\om\in \Gamma(t)$ and $h_I< \ln^M t$, we deduce that $|U(h_I)-U(h_I-1)|\leq 2K_1\ln\ln t$. 
Thus, we obtain the following upper bound for $\frac{\pi(h_I)}{\pi(a)}$,
\begin{equation}
\label{WATQ}
\frac{\pi(h_I)}{\pi(a)}\leq e^{-U(h_I)+U(a)}(2K_1+1)\ln t.
\end{equation}

\section{Technical lemmas}
\label{Technical}
We start by showing four lemmas on the asymptotic behavior of the potential $V$. 
We mention that since $V$ is defined on $\R^+$, all the intervals considered in this section are intervals of $\R^+$. Let us recall that $s(t)=(C^*(\ln\ln t)^{-1}\ln t)^{\frac{1}{\alf}}$.

\begin{figure}[!htb]
\begin{center}
\includegraphics[scale= 0.6]{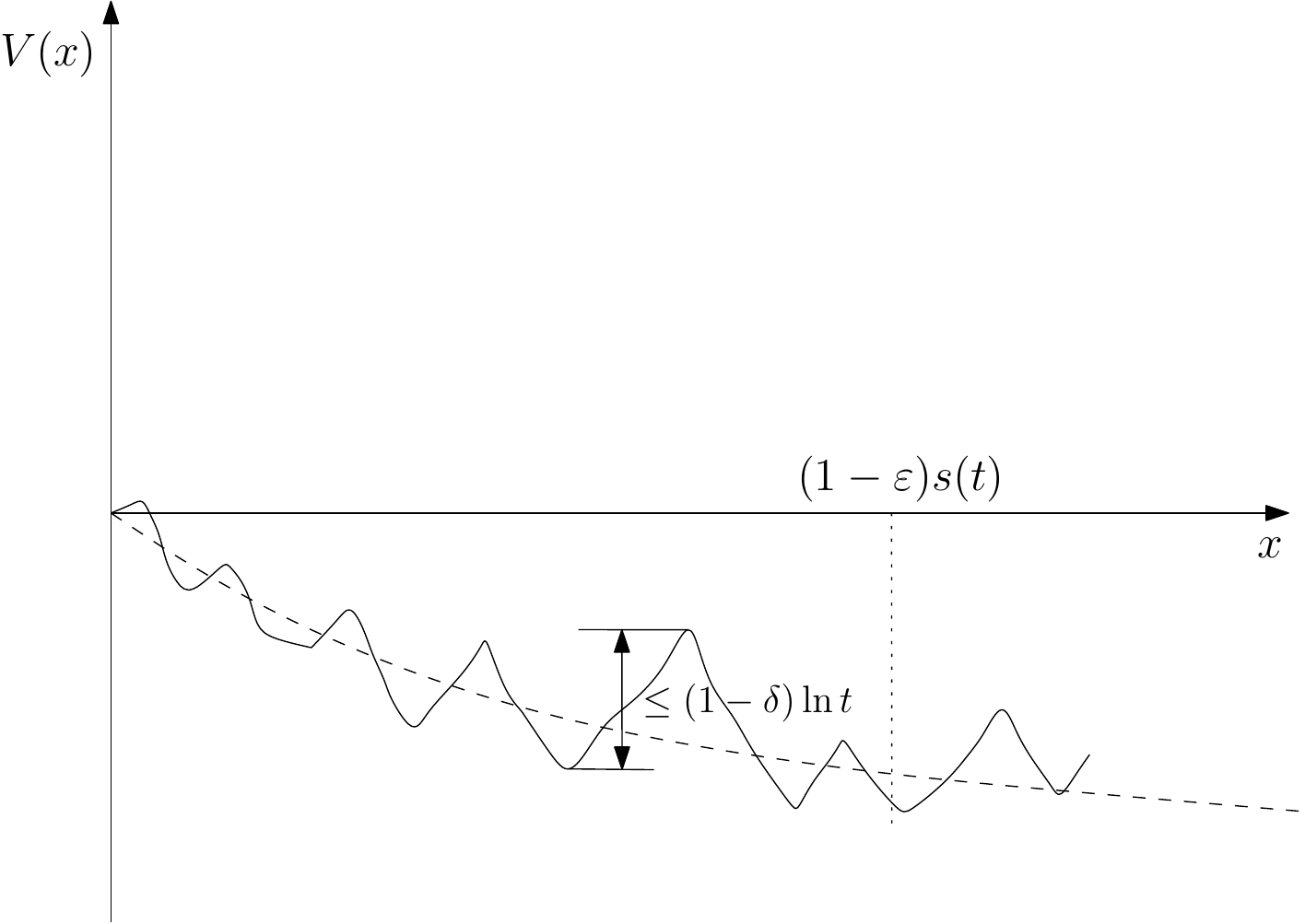}
\caption{Maximum drawup  of $V$ before $(1-\eps) s(t)$.}
\label{fig2}
\end{center}
\end{figure}

\begin{figure}[!htb]
\begin{center}
\includegraphics[scale= 0.6]{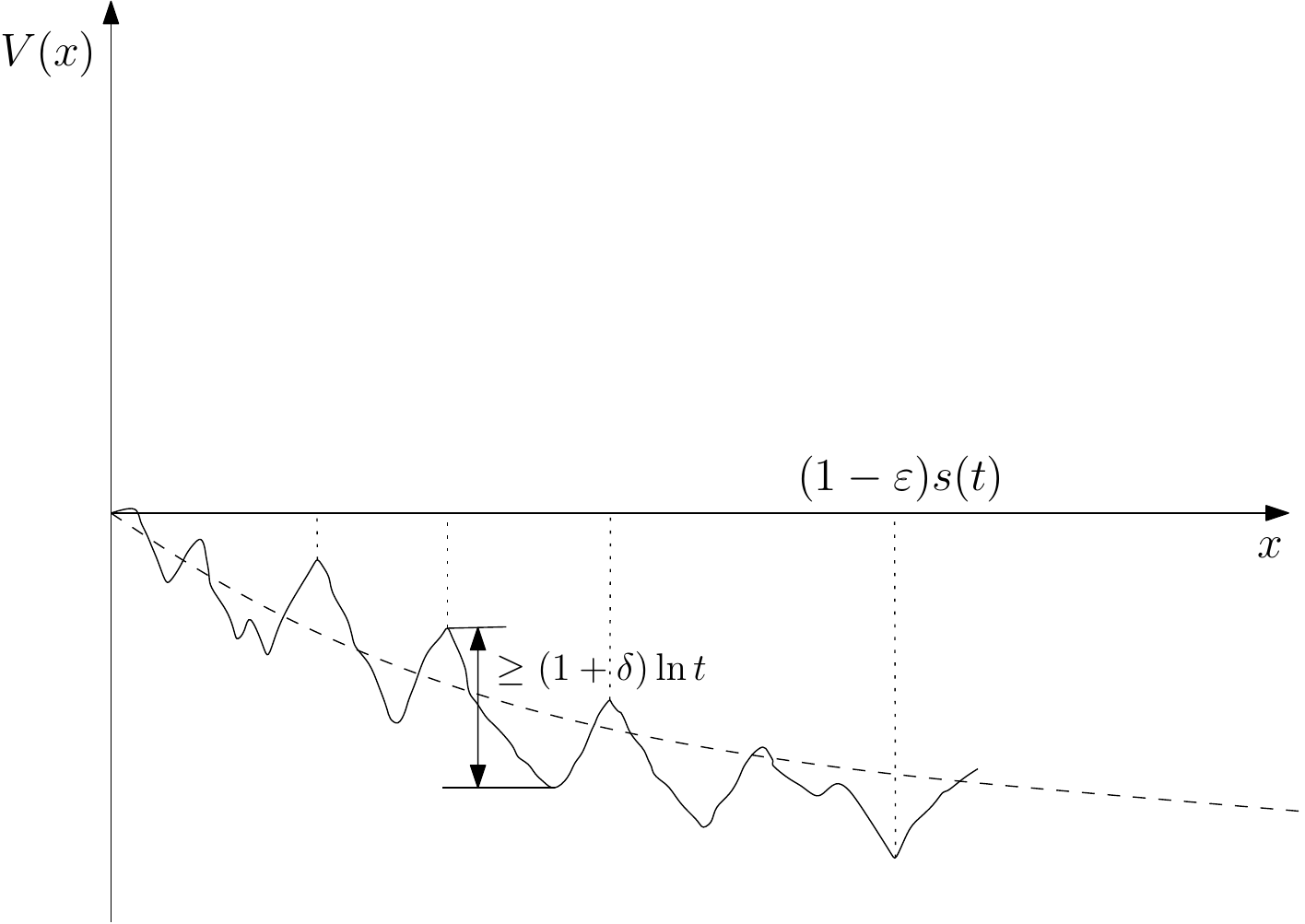}
\caption{Partition of $[0,(1-\eps) s(t)]$ into $N=4$ intervals.}
\label{fig2b}
\end{center}
\end{figure}

\begin{figure}[!htb]
\begin{center}
\includegraphics[scale= 0.6]{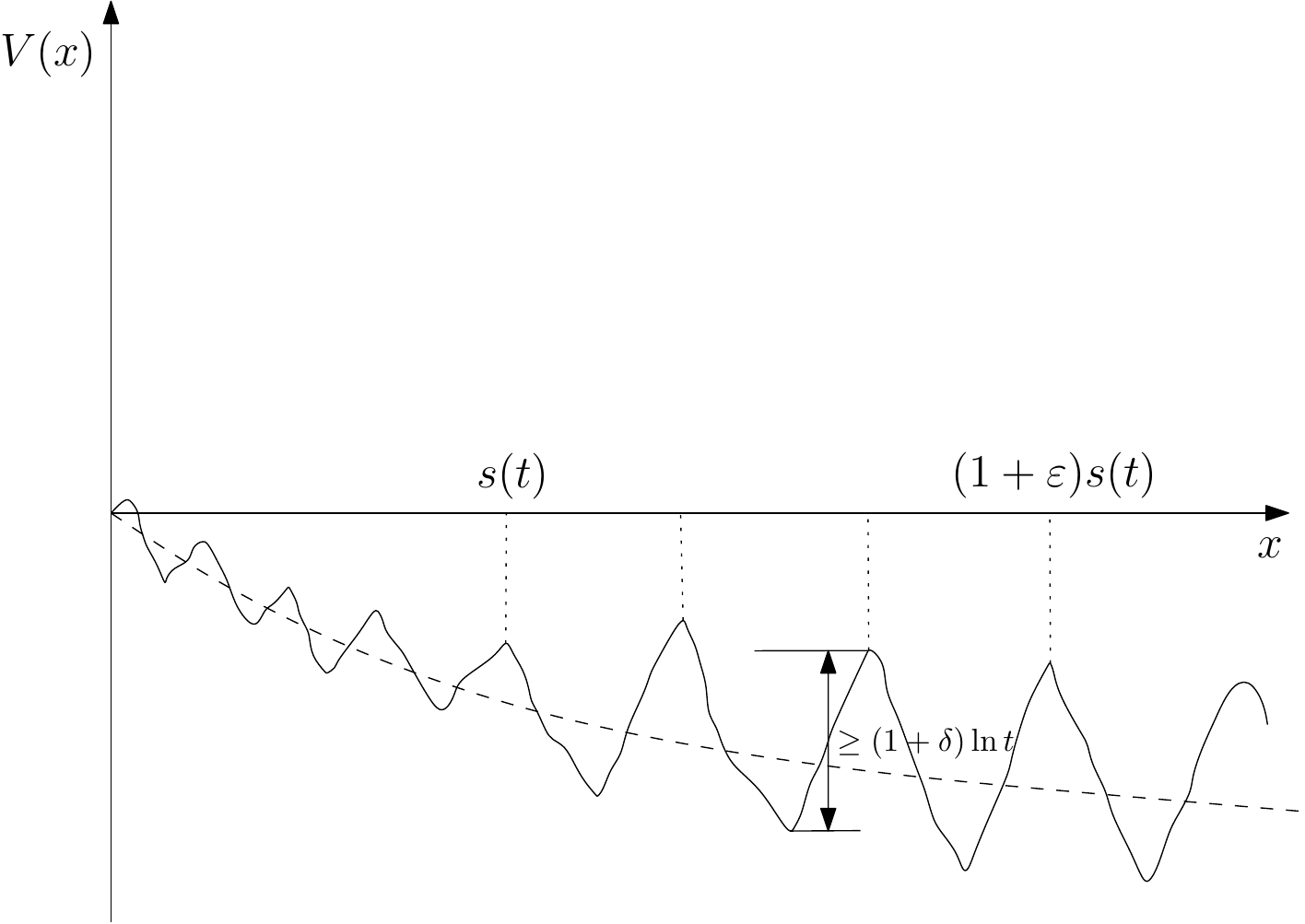}
\caption{Partition of $[s(t), (1+\eps) s(t)]$ into $N=3$ intervals.}
\label{fig2c}
\end{center}
\end{figure}

In Lemma \ref{Lem1}, we show that $\IP$-a.s., for all $t$ large enough the maximum draw-up of $V$ before $(1-\eps)s(t)$ is smaller than $(1-\delta)\ln t$, for $\delta$ suitably chosen 
(see Fig.\ \ref{fig2}). In Lemma \ref{Lem2}, we show that for any integer $N$, 
we have that, $\IP$-a.s., for all $t$ large enough, there exists a partition 
of $[0, (1-\eps)s(t)]$ into $N$ intervals such that on each interval the maximum draw-down of $V$ is greater than $(1+\delta)\ln t$ (see Fig.\ \ref{fig2b}). In Lemma \ref{Lem3}, we show that for any 
integer $N$, we have that, $\IP$-a.s., for all $t$ large enough, there exists a partition of $[s(t), (1+\eps)s(t)]$ into $N$ intervals such that on each interval the maximum draw-up of $V$ is greater than $(1+\delta)\ln t$ for $\delta$ suitably chosen (see Fig.\ \ref{fig2c}). 
Finally, in Lemma~\ref{Lem4}, we show that on the interval $[0,\ln^{\frac{1}{\alf}} t]$ the range of $V$ is smaller than $2\ln^{\frac{1}{\alf}}t$.
The proofs of Lemmas \ref{Lem2} and \ref{Lem4} follow from standard properties of 
Brownian motion. To prove Lemmas~\ref{Lem1} and \ref{Lem3} we essentially use the the same method, that is, we first approximate the potential $V$ by some suitable 
drifted Brownian motion and then apply Corollary \ref{Corro1}.

For $\eps\in (0,1)$, $\delta\in (0,1)$ and $N\in \N$, let us define the following events
\[
A_{\eps,\delta}(t):=\Big\{D^+_{[0,(1-\eps)s(t)]}(V)\leq(1-\delta)\ln t\Big\},
\] 
\begin{align*}
B_{\eps,\delta,N}(t)&:=\Big\{\mbox{there exists a partition of $[(1-\eps)s(t),(1-\frac{\eps}{2})s(t)]$ into $N$ intervals $I_j$ such that} \nonumber\\
                                  & \phantom{****}D^-_{I_j}(V)>(1+\delta)\ln t, j=1,\dots,N\Big\},
\end{align*}
and
\begin{align*}
C_{\eps,\delta,N}(t)&:=\Big\{\mbox{there exists a partition of $[s(t),(1+\eps)s(t)]$ into $N$ intervals $J_j$ such that} \nonumber\\
& \phantom{****}D^+_{J_j}(V)>(1+\delta)\ln t, j=1,\dots,N\Big\}.
\end{align*}

\begin{figure}[!htb]
\begin{center}
\includegraphics[scale= 0.6]{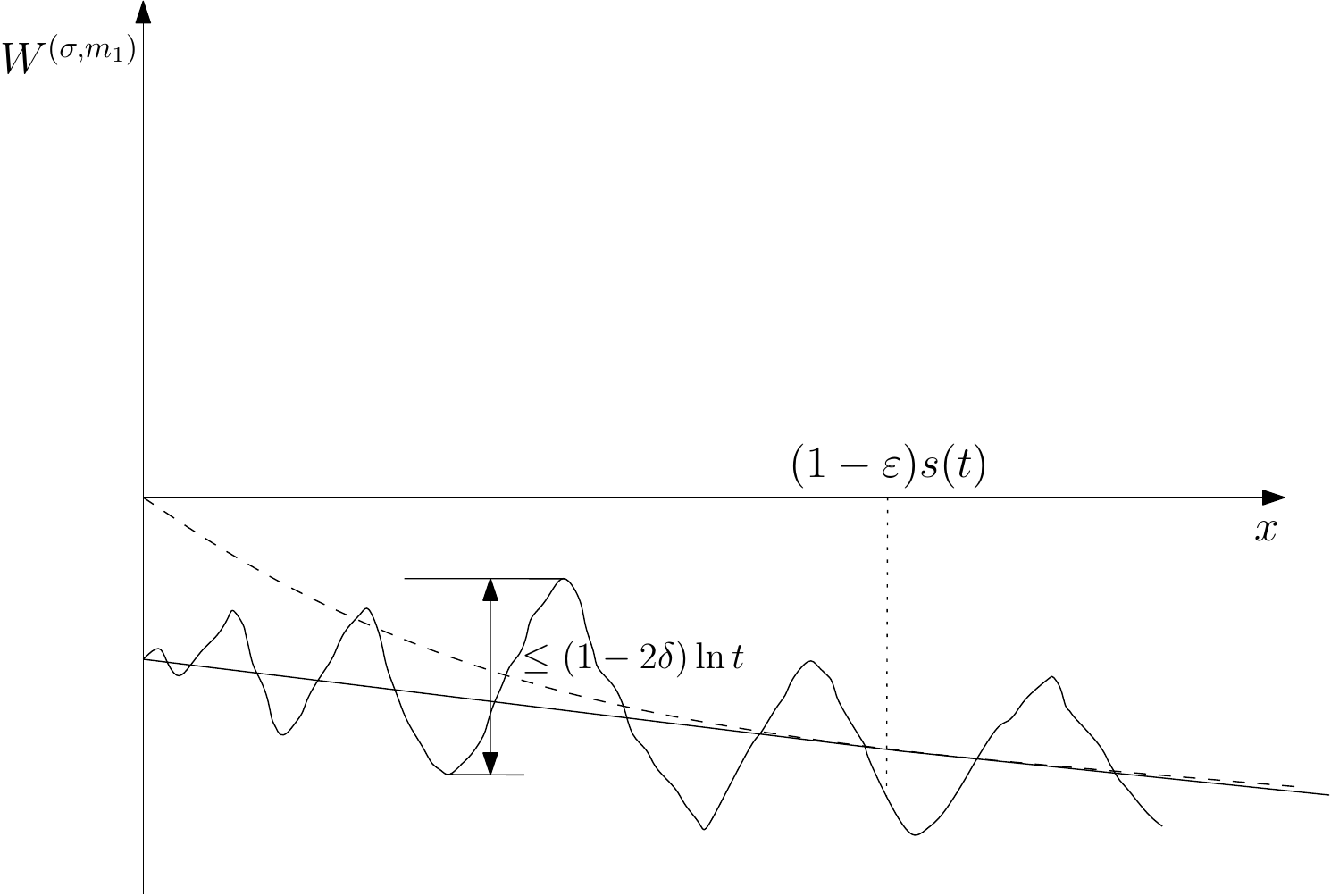}
\caption{On the definition of $W^{(\sig, m_1)}$.}
\label{fig3}
\end{center}
\end{figure}

We first show the following 
\begin{lm}
\label{Lem1}
For all $\eps\in (0,1)$, there exists $\delta>0$ small enough such that $\IP[\liminf_{t\to \infty}A_{\eps,\delta}(t)]=1$.
\end{lm}
\textit{Proof.} Consider an exponential random variable $T$ with parameter 1 and independent of $W$. Let us also introduce the drifted Brownian motion $W^{(\sig, m_1)}(x):=\sig W(x)+m_1x$ where $m_1:=-\frac{b}{(1-\eps)^{\alf}s^{\alf}(t)}$ is the derivative of the function $-\frac{b}{1-\alf}x^{1-\alf}$ at point $(1-\eps)s(t)$ (see Fig.~\ref{fig3}).
By the choice of $m_1$, we have that the event $\{D^+_{[0,(1-\eps)s(t)]}(V)>(1-2\delta)\ln t\}$ is contained in the event $\{D^+_{[0,(1-\eps)s(t)]}(W^{(\sig, m_1)})>(1-2\delta)\ln t\}$, this implies that
\begin{align}
\label{EVEnt1}
\IP[A^c_{\eps,2\delta}(t)]
&\leq \IP\Big[D^+_{[0,((1-\eps)\vee (T(\ln\ln t)^2 ))s(t)]}(W^{(\sig,m_1)})>(1-2\delta)\ln t\Big]\nonumber\\
&\leq   \IP\Big[D^+_{[0,T(\ln\ln t)^2 s(t)]}\Big(W^{(\sig,m_1)} \Big)>(1-2\delta)\ln t\Big]+\IP\Big[T\leq \frac{1-\eps}{(\ln\ln t)^2}\Big].
\end{align} 

As $T$ is exponentially distributed with parameter 1, 
the second term of the right-hand side of (\ref{EVEnt1}) is equal to 
\begin{equation}
\label{RTY1}
\IP\Big[T\leq \frac{1-\eps}{(\ln\ln t)^2}\Big]=\frac{1-\eps}{(\ln\ln t)^{2}}
\end{equation}
as $t \to \infty$.
For the first term, by Corollary \ref{Corro1} we obtain 
\begin{equation}
\label{RTY2}
 \IP\Big[D^+_{[0,Ts(t)\ln\ln t ]}\Big(W^{(\sig, m_1)} \Big)>(1-2\delta)\ln t\Big]=\frac{(1+o(1))}{1+(\ln t)^{(\frac{1}{\alf}-2)\Big(\frac{1-2\delta}{(1-\eps)^{\alf}}-1\Big)+o(1)}}
\end{equation}
as $t\to \infty$.
Now, let $\mu>0$ and consider the sequence of time intervals $I_n:=[t_n, t_{n+1})$, where $t_n:=e^{(1+\mu)^n}$ for $n\geq 0$.
Choosing $0<2\delta<1-(1-\eps)^{\alf}$ and
 using (\ref{EVEnt1}), (\ref{RTY1}) and (\ref{RTY2}) we obtain that $\sum_{n\geq 0}\IP[A^c_{\eps,2\delta}(t_n)]<\infty$.
Thus, by Borel-Cantelli Lemma we obtain that for $\IP$-a.a.\ $\om$ 
there exists $n_0=n_0(\om)$ such that $\omega \in A_{\eps,2\delta}(t_n)$ 
for all $n\geq n_0$. Now, let $n\geq n_0$ and 
suppose $t\in [t_n, t_{n+1})$. We have $\IP$-a.s.,
\begin{align*}
D^+_{[0,(1-\eps)s(t)]}(V)&\leq D^+_{[0,(1-\eps)s(t_{n+1})]}(V)\nonumber\\
&\leq (1-2\delta)\ln t_{n+1}\nonumber\\
&=(1-2\delta)(1+\mu)\ln t_n\nonumber\\
&\leq (1-2\delta)(1+\mu)\ln t.
\end{align*}
Choosing $\mu$ in such a way that $(1-2\delta)(1+\mu)\leq (1-\delta)$, 
we obtain that for $\IP$-a.a.\ $\om$ there exists $t_0=t_0(\om)$ such that 
$\omega \in A_{\eps,\delta}(t)$ for all $t\geq t_0$, which proves Lemma \ref{Lem1}.
\qed

\begin{lm}
\label{Lem2}
For all $\eps\in (0,1)$ and  $\delta>0$ we have $\IP[\liminf_{t\to \infty}B_{\eps,\delta,N}(t)]=1$, for all $N\geq 1$.
\end{lm}
\textit{Proof.} Let $\mu>0$ be such that $\beta:=(1-\eps)(1+\mu)^{\frac{1}{\alf}}< (1-\frac{\eps}{2})$ and consider the sequence of time intervals $I_n:=[t_n, t_{n+1})$, 
where $t_n:=e^{(1+\mu)^n}$ for $n\geq 0$. Divide the interval $[\beta s(t),(1-\frac{\eps}{2})s(t)]$ into $N$ intervals $\mathcal{I}_j$, $j=1,\dots,N$, 
of size $\eta s(t)$ with $\eta:=N^{-1}(1-\frac{\eps}{2}-\beta)$. Let us define the following events
\[
E_{\eps,\delta,\mu}(t):=\bigcup_{j=1}^N\Big\{D^-_{\mathcal{I}_j}(V)\leq (1+\delta)\ln t\Big\}.
\]
We have 
\begin{align}
\label{WER}
\IP[E_{\eps,2\delta,\mu}(t)]&\leq \sum_{j=1}^{N}\IP\Big[D^-_{\mathcal{I}_j}(\sig W)
\leq (1+2\delta)\ln t\Big]\nonumber\\
&= N\IP\Big[\max_{s\in [0,\eta s(t)]}|W(s)|\leq \frac{1+2\delta}{\sig}\ln t\Big]\nonumber\\
&\leq N\IP\Big[\max_{s\in [0,\eta s(t)]}W(s)\leq \frac{1+2\delta}{\sig}\ln t\Big]\nonumber\\
&=N\Big(1-2\IP\Big[W(\eta s(t))> \frac{1+2\delta}{\sig}\ln t\Big]\Big)\nonumber\\
&=N\Big(1-2\int_{\frac{(1+2\delta)\ln t}{\sig (\eta s(t))^{1/2}}}^{\infty} \frac{e^{-\frac{y^2}{2}}}{\sqrt{2\pi}}dy\Big)\nonumber\\
&=N\sqrt{\frac{2}{\pi}}\frac{1+2\delta}{\sig\eta^{\frac{1}{2}} (C^*)^{\frac{1}{2\alf}}}(\ln\ln t)^{\frac{1}{2\alf}}(\ln t)^{-(\frac{1}{2\alf}-1)} (1+o(1))
\end{align}
as $t\to \infty$.
We obtain from (\ref{WER}) that $\sum_{n\geq 0}\IP[E_{\eps,2\delta,\mu}(t_n)]<\infty$. 
Thus, by Borel-Cantelli Lemma we obtain that for $\IP$-a.a.\ $\om$ there exists $n_0=n_0(\om)$ such that $\omega \in E^c_{\eps,2\delta,\mu}(t_n)$ for all $n\geq n_0$. 
Now, suppose that $n\geq n_0$ and $t\in [t_n, t_{n+1})$. 
Since we have $s^{\alf}(t_n)\leq s^{\alf}(t)\leq (1+\mu)s^{\alf}(t_n)$ for large enough $n$, we deduce that $\IP$-a.s., 
there exists a partition of $[(1-\eps)s(t),(1-\frac{\eps}{2})s(t)]$ into $N$ intervals $I_j$, $j=1,\dots,N$, such that on each one $D^-_{I_j}(V)> (1+2\delta)\ln t_n$.
%\begin{equation}
%\label{WERT}
%D^-_{[(1-\eps)s(t),(1-\frac{\eps}{2})s(t)]}(V)\geq D^-_{[\beta s(t_n),(1-\frac{\eps}{2})s(t_n)]}(V)>(1+2\delta)\ln t_n.
%\end{equation}
Since $\ln t_n\leq \ln t\leq (1+\mu)\ln t_n$, 
we have $(1+2\delta)\ln t_n\geq \frac{1+2\delta}{1+\mu}\ln t\geq (1+\delta)\ln t$ for $\mu>0$ 
small enough. From these last observations, we conclude that for $\IP$-a.a.\ $\om$, there exists $t_0=t_0(\om)$ such that $\omega \in B_{\eps,\delta,N}(t)$ for all $t\geq t_0$, which proves Lemma \ref{Lem2}.
\qed

\begin{figure}[!htb]
\begin{center}
\includegraphics[scale= 0.7]{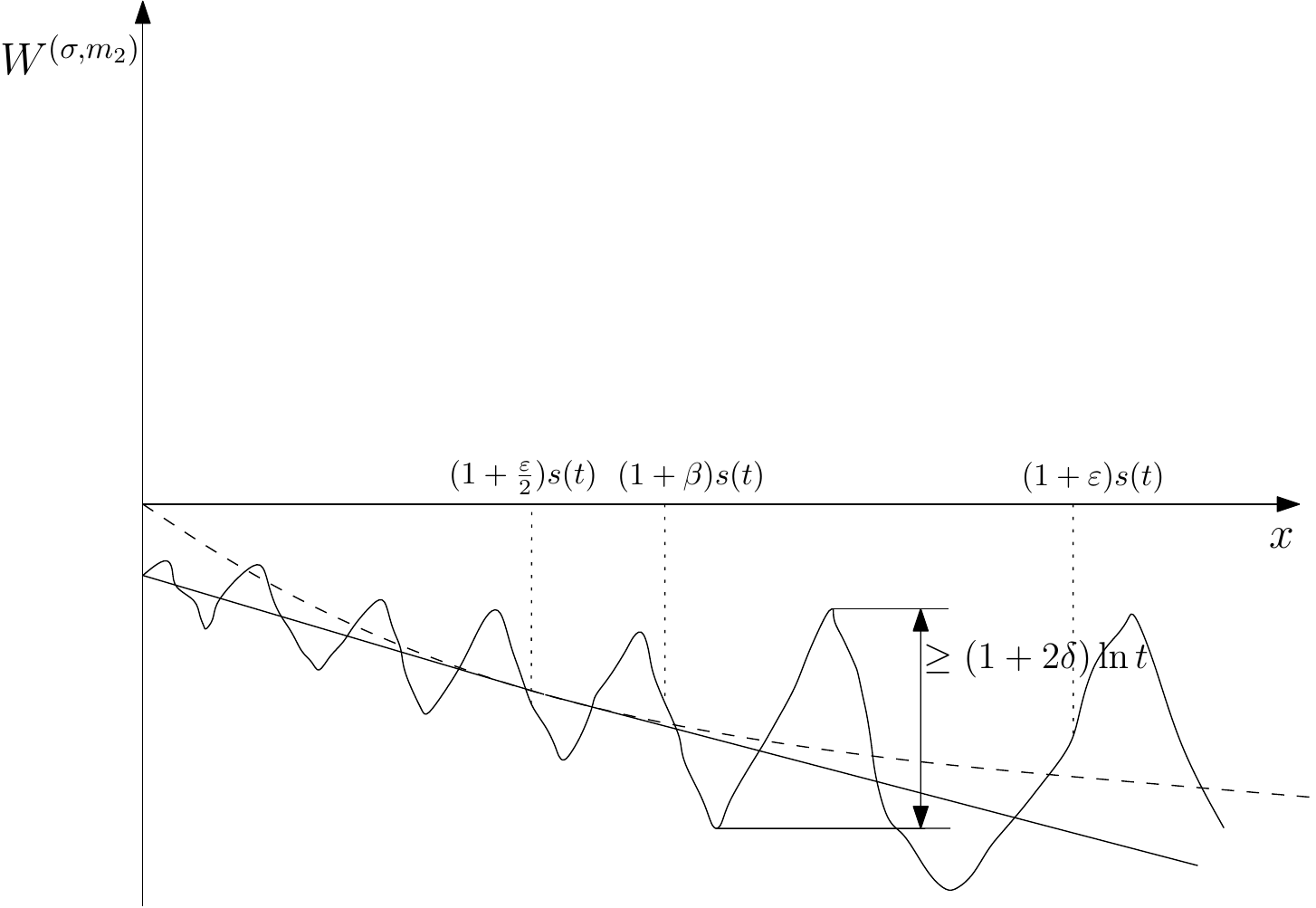}
\caption{On the definition of $W^{(\sig, m_2)}$.}
\label{fig4}
\end{center}
\end{figure}

\begin{lm}
\label{Lem3}
For all $\eps\in (0,1)$, there exists small enough $\delta>0$ such that 
$\IP[\liminf_{t\to \infty}C_{\eps,\delta,N}(t)]=1$, for all $N\geq 1$.
\end{lm}
\textit{Proof.} Let $\mu>0$ be such that $(1+\beta):=(1+\frac{\eps}{2})(1+\mu)^{\frac{1}{\alf}}< (1+\eps)$ and consider again the sequence of time intervals $I_n:=[t_n, t_{n+1})$, 
where $t_n=e^{(1+\mu)^n}$ for $n\geq 0$. Divide the interval $[(1+\beta) s(t),(1+\eps)s(t)]$ into $N$ intervals $\mathcal{J}_j$, $j=1,\dots,N$ of size $\frac{\eps-\beta}{N} s(t)$. Let us define the following events
\[
F_{\eps,\delta,\mu}(t):=\bigcup_{j=1}^N
\Big\{D^+_{\mathcal{J}_j}(V)\leq(1+\delta)\ln t\Big\}.
\]
Let $m_2:=-\frac{b}{(1+2^{-1}\eps)^{\alf}s^{\alf}(t)}$ 
be the derivative of the function $-\frac{b}{1-\alf}x^{1-\alf}$ at 
point $(1-\frac{\eps}{2})s(t)$ and introduce the drifted Brownian motion $W^{(\sig, m_2)}(x):=\sig W(x)+m_2 x$ (see Fig.~\ref{fig4}). By definition of $W^{(\sig, m_2)}$, 
we have that the event $\Big\{D^+_{[(1+\beta)s(t),(1+\beta+T(\ln\ln t)^{-1})s(t)]}(V)\leq(1+2\delta)\ln t\Big\}$ is contained in the event 
$\Big\{D^+_{[(1+\beta)s(t),(1+\beta+T(\ln\ln t)^{-1})s(t)]}(W^{(\sig, m_2)})\leq(1+2\delta)\ln t\Big\}$, this leads to
\begin{align}
\label{event5}
\IP[F_{\eps,2\delta,\mu}(t)]&\leq \sum_{j=1}^{N} 
\IP\Big[D^+_{\mathcal{J}_j}(V)\leq(1+2\delta)\ln t\Big]\nonumber\\
&\leq N\IP\Big[D^+_{[(1+\beta)s(t),((1+\beta+N^{-1}(\eps-\beta))\wedge (1+\beta+T(\ln\ln t)^{-1}))s(t)]}(V)\leq (1+2\delta)\ln t\Big]\nonumber\\
&\leq N \IP\Big[D^+_{[(1+\beta)s(t),(1+\beta+T(\ln\ln t)^{-1})s(t)]}(V)\leq(1+2\delta)\ln t\Big]+N\IP\Big[T>\frac{\eps-\beta}{N}\ln\ln t \Big]\nonumber\\
&\leq N \IP\Big[D^+_{[(1+\beta)s(t),(1+\beta+T(\ln\ln t)^{-1})s(t)]}(W^{(\sig, m_2)})\leq(1+2\delta)\ln t\Big]+N\IP\Big[T>\frac{\eps-\beta}{N}\ln\ln t \Big]\nonumber\\
&=   N\IP\Big[D^+_{[0,Ts(t)(\ln\ln t)^{-1} ]}\Big( W^{(\sig, m_2)} \Big)\leq(1+2\delta)\ln t\Big]+N\IP\Big[T>\frac{\eps-\beta}{N}\ln\ln t \Big].
\end{align}

As  $T$ is exponentially distributed with parameter 1, 
we have for the second term of the right-hand side of (\ref{event5})
\begin{equation}
\label{event6}
N\IP\Big[T>\frac{\eps-\beta}{N}\ln\ln t \Big]=N\ln^{\frac{\eps-\beta}{N}} t. 
\end{equation}
For the first term, we use Corollary \ref{Corro1} to obtain that
\begin{equation}
\label{event7}
\IP\Big[D^+_{[0,Ts(t)(\ln\ln t)^{-1} ]}\Big(W^{(\sig, m_2)} \Big)\leq(1+2\delta)\ln t\Big]=1-\frac{(1+o(1))}{1+(\ln t)^{(\frac{1}{\alf}-2)\Big(\frac{1+2\delta}{(1+2^{-1}\eps)^{\alf}}-1\Big)+o(1)}}
\end{equation}
as $t\to \infty$.
 Choosing $0<2\delta<(1+2^{-1}\eps)^{\alf}-1$ and using (\ref{event5}), (\ref{event6}) and (\ref{event7}) we obtain that $\sum_{n\geq 0}\IP[F_{\eps,2\delta,\mu}(t_n)]<\infty$.
Thus, by 
Borel-Cantelli Lemma we obtain that for $\IP$-a.a.\ $\om$ there exists $n_0=n_0(\om)$ such that $\omega \in F^c_{\eps,2\delta,\mu}(t_n)$ for all $n\geq n_0$. Now, let $n\geq n_0$ and 
suppose $t\in [t_n, t_{n+1})$. 
Since we have $s^{\alf}(t_n)\leq s^{\alf}(t)\leq (1+\mu)s^{\alf}(t_n)$, we deduce that $\IP$-a.s., there exists a partition of $[(1+\frac{\eps}{2})s(t),(1+\eps)s(t)]$ 
into $N$ intervals $J_j$, $j=1,\dots,N$, such that on each one $D^+_{J_j}(V)> (1+2\delta)\ln t_n$. 
As $\ln t_n\leq \ln t\leq (1+\mu)\ln t_n$, we have $(1+2\delta)\ln t_n\geq \frac{1+2\delta}{1+\mu}\ln t\geq (1+\delta)\ln t$ for $\mu>0$ 
small enough. From these last observations, we conclude that for $\IP$-a.a.\ $\om$, there exists $t_0=t_0(\om)$ such that $\omega \in C_{\eps,\delta,N}(t)$ for all $t\geq t_0$, which proves Lemma \ref{Lem3}.
\qed
Finally, let $G(t):=\Big\{\max_{y\leq \ln^{1/\alf} t}|V(y)|\leq 2\ln^{\frac{1}{\alf}}t\Big\}$. We show the following

\begin{lm}
\label{Lem4}
We have that $\IP[\liminf_{t\to \infty}G(t)]=1$. 
\end{lm}
\textit{Proof.} Let $n$ be an positive integer. By \cite{PerMot}, Lemma 12.9, we have
\begin{align*}
\IP\Big[\max_{y \leq \ln^{1/ \alf} (n+1)}|V(y)|> 2\ln^{\frac{1}{\alf}}n\Big]&\leq 
\IP\Big[\max_{y\leq  \ln^{1/\alf} (n+1)}|W(y)|>\sig^{-1}\ln^{\frac{1}{\alf}}n\Big]\nonumber\\
&\leq 2\IP\Big[\max_{y\leq \ln^{1/ \alf} (n+1)}W(y)>\sig^{-1}\ln^{\frac{1}{\alf}}n\Big]\nonumber\\
&=2\IP\Big[W( \ln^{\frac{1}{\alf}} (n+1))>\sig^{-1}\ln^{\frac{1}{\alf}}n\Big]\nonumber\\
&\leq \frac{2}{\sqrt{2\pi \sig^2}}e^{-\frac{\ln^{\frac{1}{\alf}} (n+1)}{{2\sig^2}}}
\end{align*}
for sufficiently large $n$. Since $\alf\in(0,\frac{1}{2})$, we deduce that $\sum_{n>1}\IP\Big[\max_{y\leq \ln^{1/\alf} (n+1)}|V(y)|> 2\ln^{\frac{1}{\alf}}n\Big]<\infty$. 
By Borel-Cantelli Lemma, we have that for $\IP$-a.a.\ $\om$ there exists $n_0=n_0(\om)$ such that for all $n\geq n_0$ we have $\max_{y\leq \ln^{1/\alf} (n+1)}|V(y)|\leq 2\ln^{\frac{1}{\alf}}n$. 
Now consider $n\geq n_0$ and $t\in [n,n+1)$, we have that $\max_{y\leq \ln^{1/\alf} t}|V(y)|\leq \max_{y\leq \ln^{1/\alf} (n+1)}|V(y)|\leq 2\ln^{\frac{1}{\alf}}n\leq 2\ln^{\frac{1}{\alf}}t$. 
This shows that $\IP[\liminf_{t\to \infty}G(t)]=1$ and concludes the proof of Lemma \ref{Lem4}.
\qed

\section{Proof of Theorem~\ref{Theo}}
\label{secTheo}
In this last section, for the sake of brevity, expressions like $X_t=x$ or $\tau_x>t$ must be understood as $X_t=\lf x\rf$ or $\tau_{\lf x \rf}>t$ (where $\lf \cdot\rf$ is the integer part function) whenever $x$ in not necessarily integer. 
Also, in contrast with the former section, all the intervals considered in this section are intervals of $\Z^+$. We will also need the function $\lc\cdot \rc:=\lf \cdot\rf+1$.
\medskip

Fix some $\eps\in (0,1)$. We start by showing that for $\IP$-a.a.\ $\om$, $\Po[\liminf_{t \to \infty}s(t)^{-1}X_t\geq (1-\eps)]=1$. Let $\delta\in (0,1)$ be such that Lemmas \ref{Lem1} and \ref{Lem2} hold. 
Take $N=\lf 2\delta^{-1}\rf$ and let $\om$ be such that $\om\in \liminf_{t\to \infty}(A_{\eps, \delta}(t)\cap B_{\eps, \delta,N}(t)\cap G(t)\cap \Gamma(t))$.
Let us define 
\[
\hat{\tau}(t):=\inf\{u> \tau_{\lceil(1-\frac{\eps}{2})s(\lf t\rf)\rceil}: X_u=(1-\eps)s(\lc t\rc)\}
\]
for all $t\geq 3$, with the convention $\inf\{\emptyset\}=\infty$.
We have for all integer $n\geq 3$,
\begin{align}
\label{FUG}
\Po[\{\tau_{\lc(1-\frac{\eps}{2})s(n)\rc} \geq n\}\cup\{\hat{\tau}(n)-\tau_{\lc(1-\frac{\eps}{2})s(n)\rc}\leq n\}] &\leq \Po[\tau_{\lc(1-\frac{\eps}{2})s(n)\rc} \geq n]\nonumber\\
&\phantom{**}+\Po[\hat{\tau}(n)-\tau_{\lc(1-\frac{\eps}{2})s(n)\rc}\leq n].
\end{align}
The next step is to apply Proposition \ref{Confine} 
to the first term of the right-hand side of (\ref{FUG}). Since $\om\in \liminf_{t\to \infty}(A_{\eps, \delta}(t)\cap G(t)\cap \Gamma(t))$, 
we have that for $n$ large enough $H([0,\lc(1-\frac{\eps}{2})s(n)\rc]\leq D^+_{[0,\lc(1-\frac{\eps}{2})s(n)\rc]}(U)\leq (1-\delta)\ln n+o(\ln n)$ and  $\tilde{M}\leq 4\ln^{\frac{1}{\alf}} n+o(\ln n)$. 
Therefore, by Proposition~\ref{Confine}
we obtain
\begin{align}
\label{FUG1}
 \Po[ \tau_{\lc(1-\frac{\eps}{2})s(n)\rc} \geq n]\leq \exp{\{-n^{\delta+o(1)}\}}
\end{align}
as $n \to \infty$.
For the second term of the right-hand side of (\ref{FUG}), we have by the Markov property applied at time
$\tau_{\lc(1-\frac{\eps}{2})s(n)\rc}$,
\begin{equation}
\Po[\hat{\tau}(n)-\tau_{\lc(1-\frac{\eps}{2})s(n)\rc}\leq n]= \Po^{\lc(1-\frac{\eps}{2})s(n)\rc}[\tau_{(1-\eps)s(n+1)}\leq n].
\end{equation}
Since $\om\in \liminf_{t\to \infty}B_{\eps,\delta,N}(t)\cap \Gamma(t)$ 
there exists for $n$ large enough a partition 
$x_0=\lf(1-\eps)s(n+1)\rf<x_1<\dots<x_{N-1}<x_N=\lc(1-\frac{\eps} {2})s(n)\rc$ of $[\lf(1-\eps)s(n+1)\rf,\lc(1-\frac{\eps}{2})s(n)\rc]$ into $N=\lf 2\delta^{-1}\rf$ intervals $I_j=[x_{j-1},x_{j}]$, $j=1,\dots,N$, such that on each interval $D^-_{I_j}(U)>(1+\delta)\ln n-o(\ln n)$. 
By the Markov property we have
\begin{align*}
\Po^{\lc(1-\frac{\eps}{2})s(n)\rc}[\tau_{(1-\eps)s(n+1)}\leq n]
&\leq \Po^{\lc(1-\frac{\eps}{2})s(n)\rc}[\tau_{x_{j-1}}\leq n, j=1,\dots,N]\nonumber\\
&\leq \prod_{j=1}^{N}\Po^{x_{j}}[\tau_{x_{j-1}}\leq n].
\end{align*}
Applying Proposition \ref{escape} to the right-hand side of the last inequality and using bound (\ref{WATQ}), we obtain 
\begin{align}
\label{FUG2}
\Po^{\lc(1-\frac{\eps}{2})s(n)\rc}[\tau_{(1-\eps)s(n+1)}\leq n]&\leq K_3^{\lf 2\delta^{-1}\rf}(n+1)^{-(2-\delta)+o(1)}
\end{align}
as $n\to \infty$.
From (\ref{FUG}), (\ref{FUG1}) and (\ref{FUG2}), 
as $\delta\in (0,1)$, we deduce that $\sum_{n\geq 3}\Po[\{\tau_{(1-\frac{\eps}{2})s(n)} \geq n\}\cup\{\hat{\tau}(n)-\tau_{(1-\frac{\eps}{2})s(n)}\leq n\}]<\infty$. 
By Borel-Cantelli Lemma, we obtain that, $\Po$-a.s., for all $n$ large enough $X_n>(1-\eps)s(n)$. 
Now, for $t\in [n,n+1)$ and $n$ large enough, we have that $\tau_{\lc(1-\frac{\eps}{2})s(n)\rc}<n\leq t$ and $\hat{\tau}(t)-\tau_{\lc(1-\frac{\eps}{2})s(n)\rc}\geq n+1>t$, 
which implies $X_t> (1-\eps)s(t)$. By Lemmas \ref{Lem1}, \ref{Lem2}, \ref{Lem4} 
and the definition of $\Gamma(t)$, we conclude that for $\IP$-a.a.\ $\om$, $\Po[\liminf_{t \to \infty} s(t)^{-1}X_t\geq (1-\eps)]=1$.
\medskip

We continue the proof of Theorem \ref{Theo} by showing that for $\IP$-a.a.\ $\om$, $\Po[\limsup_{t \to \infty} s(t)^{-1}X_t\leq(1+\eps)]=1$.
Let $\delta\in (0,1)$ be such that Lemma \ref{Lem3} holds, $N=\lf 2\delta^{-1}\rf$ and $\om$ be such that $\om\in \liminf_{t\to \infty}(C_{\eps, \delta,N}(t)\cap \Gamma(t))$. 
Since $\om\in \liminf_{t\to \infty}(C_{\eps, \delta,N}(t)\cap \Gamma(t))$ 
there exists for all large enough integers $n$ a partition $y_0=0<y_1<\dots<y_{N-1}<y_N=\lf(1+\eps)s(n)\rf$ of $[0,\lf(1+\eps)s(n)\rf]$ into $N=\lf 2\delta^{-1}\rf$ intervals $J_j=[y_{j-1},y_{j}]$, $j=1,\dots,N$, 
such that on each interval $D^+_{J_j}(U)>(1+\delta)\ln n-o(\ln n)$. By the Markov property we have
\begin{align*}
\Po[\tau_{(1+\eps)s(n)}\leq n]\leq \prod_{j=1}^{N}\Po^{y_{j-1}}[\tau_{y_{j}}\leq n].
\end{align*}
Applying Proposition \ref{escape} to the right-hand term of the last inequality and using bound (\ref{WATQ}), we obtain
\begin{align}
\label{ERT}
\Po[\tau_{(1+\eps)s(n)}\leq n]\leq K_3^{\lf 2\delta^{-1}\rf}(n+1)^{-(2-\delta)+o(1)}
\end{align}
as $n\to \infty$.
From (\ref{ERT}), as $\delta\in (0,1)$, we deduce that $\sum_{n\geq 3}\Po[\tau_{(1+\eps)s(n)}\leq n]<\infty$. 
By Borel-Cantelli Lemma, we obtain that, $\Po$-a.s., for all $n$ large enough $X_n<(1+\eps)s(n)$. 
Now, for $t\in [n,n+1)$ and $n$ large enough, we have that $\tau_{(1+\eps)s(t)}\geq \tau_{(1+\eps)s(n)}\geq n+1> t$, which implies $X_t<(1+\eps)s(t)$. 
By Lemma~\ref{Lem3} and the definition of $\Gamma(t)$, we conclude that for $\IP$-a.a.\ $\om$, 
$$\Po\Big[\limsup_{t \to \infty}\frac{X_t}{s(t)}\leq(1+\eps)\Big]=1.$$
\medskip

To sum up, we showed that for $\IP$-a.a.\ $\om$, 
\[
\Po\Big[\liminf_{t\to \infty}\frac{X_t}{s(t)}\geq(1-\eps), 
\limsup_{t \to \infty}\frac{X_t}{s(t)}\leq(1+\eps)\Big]=1.
\]
 As $\eps$ is arbitrary, this shows Theorem \ref{Theo}.
\qed

%%%%%%%%%%%%%%%%%%%%%%%%%%%%%%%%%%%%%%%%%%%%%%%%%%%%%%%%%%%%%%%%%% 

\section*{Acknowledgements}
C.G.\ is grateful to FAPESP (grant 2009/51139--3) for financial
support. G.M.S.\ thanks FAPESP (grant 2011/21089-4) and S.P.\  thanks 
CNPq (grant 301644/2011-0) for financial supports. C.G.\ and S.P.\  thank 
FAPESP (grant 2009/52379-8) for financial support. 
C.G.\ and G.M.S.\ thank NUMEC for kind hospitality.

\end{document}